\input amstex
\documentstyle{amsppt}
\magnification=\magstep 1
\loadbold

\pagewidth{6.4truein}
\pageheight{8.5truein}

\define\A{{\Cal A}}
\redefine\B{{\Bbb B}}
\define\C{{\Bbb C}}
\redefine\D{{\Bbb D}}
\redefine\O{{\Cal O}}
\redefine\P{{\Bbb P}}
\define\R{{\Bbb R}}
\redefine\S{{\Bbb S}}

\define\g{{\gamma}}
\define\G{{\Gamma}}
\redefine\epsilon{{\varepsilon}}
\redefine\d{{\partial}}

\define\Comm{\operatorname{Comm}}
\redefine\Im{\operatorname{Im}}
\redefine\Re{\operatorname{Re}}
\define\Ric{\operatorname{Ric}}
\define\ac{\acuteaccent}
\define\gr{\graveaccent}
\redefine\cdot{\boldsymbol\cdot}

\NoRunningHeads
\refstyle{A} 

\topmatter 
\title Holomorphic functions of slow growth \\ 
       on nested covering spaces of compact manifolds   \endtitle
\author Finnur L\ac arusson \endauthor 
\affil University of Western Ontario  \endaffil 
\address Department of Mathematics, University of Western Ontario,
	London, Ontario N6A 5B7, Canada \endaddress
\email larusson\@uwo.ca \endemail 
\date 3 March 1999  \enddate 

\thanks This work was supported in part by the Natural Sciences and
Engineering Research Council of Canada.  \endthanks


\abstract Let $Y$ be an infinite covering space of a projective manifold
$M$ in $\P^N$ of dimension $n\geq 2$.  Let $C$ be the intersection with
$M$ of at most $n-1$ generic hypersurfaces of degree $d$ in $\Bbb P^N$. 
The preimage $X$ of $C$ in $Y$ is a connected submanifold.  Let $\phi$
be the smoothed distance from a fixed point in $Y$ in a metric pulled up
from $M$.  Let $\O_\phi(X)$ be the Hilbert space of holomorphic
functions $f$ on $X$ such that $f^2 e^{-\phi}$ is integrable on $X$, and
define $\O_\phi(Y)$ similarly.  Our main result is that (under more
general hypotheses than described here) the restriction $\O_\phi(Y)
\to \O_\phi(X)$ is an isomorphism for $d$ large enough. 

This yields new examples of Riemann surfaces and domains of holomorphy
in $\C^n$ with corona.  We consider the important special case when $Y$
is the unit ball $\B$ in $\C^n$, and show that for $d$ large enough,
every bounded holomorphic function on $X$ extends to a unique function
in the intersection of all the nontrivial weighted Bergman spaces on
$\B$.  Finally, assuming that the covering group is arithmetic, we
establish three dichotomies concerning the extension of bounded
holomorphic and harmonic functions from $X$ to $\B$.  \endabstract

\endtopmatter

\document

\specialhead Introduction  \endspecialhead

\flushpar Let $Y\to M$ be an infinite covering space of an
$n$-dimensional projective manifold, $n\geq 2$.  The function theory of
such spaces is still not well understood.  The central problem in this
area is the conjecture of Shafarevich that the universal covering space
of any projective manifold is holomorphically convex.  This is a
higher-dimensional variation on the venerable theme of uniformization. 
There are no known counterexamples to the conjecture, and it has been
verified only in a number of fairly special cases. 

Suppose $M$ is embedded into a projective space by sections of a very
ample line bundle $L$.  The generic linear subspace of codimension $n-1$
intersects $M$ in a 1-dimensional connected submanifold $C$ called an
$L$-curve.  The preimage $X$ of $C$ is a connected Riemann surface
embedded in $Y$.  A natural approach to constructing holomorphic
functions on $Y$ is to extend them from $X$.  This has the advantage of
reducing certain questions to the 1-dimensional case, but the price one
pays is having to work with functions of slow growth.  Here, slow growth
means slow exponential growth with respect to the distance from a fixed
base point or a similar well-behaved exhaustion, in an $L^2$ or
$L^\infty$ sense.  Functions in the Hardy class $H^p(X)$ grow slowly in
this sense for $p$ large enough. 

In Section 1, we improve upon the main result of our earlier paper
\cite{L\ac ar1} and show that if $L$ is sufficiently ample, then the
restriction map $\cdot|X$ is an isomorphism of the Hilbert spaces of
holomorphic functions of slow growth.  As before, the proof is based on
the $L^2$ method of solving the $\bar\d$-equation.  This may be viewed
as a sampling and interpolation theorem, related to those of Seip;
Berndtsson and Ortega Cerd\gr a; and others.  See \cite{BO}, \cite{Sei},
and the references therein. 

In Section 2, we use the isomorphism theorem to construct new examples
of Riemann surfaces with corona.  These easily defined surfaces have
many symmetries, we have a simple description of characters in the
corona, and the corona is large in the sense that it contains a domain
in euclidean space of arbitrarily high dimension. 

In Section 3, we adapt results of H\"ormander on generating algebras of
holomorphic functions of exponential growth to the case of covering
spaces.  Consequently, as shown in Section 4, the restriction map
$\cdot|X$ may fail to be an isomorphism if $L$ is not sufficiently
ample compared to the exhaustion. 

Under the mild assumption that the covering group is Gromov hyperbolic,
we found in \cite{L\ac ar2} that the only obstruction to every positive
harmonic function on $X$ being the real part of a holomorphic function
(in which case $X$ has many holomorphic functions of slow growth that
extend to $Y$) is a geometric condition involving the Martin boundary,
characteristic of the higher-dimensional case.  There are examples of
infinitely connected $X$ for which the obstruction is not present, but
these have 1-dimensional boundary, whereas in general the curves $X$ of
interest to us do not: they have the same boundary as the ambient space
$Y$.  No examples with higher dimensional boundary are known. 

In hopes of shedding some light on the dichotomy in \cite{L\ac ar2}, we
restrict ourselves from Section 4 onwards to what seems to be the most
auspicious setting possible and let $Y$ be the unit ball $\B$ in $\C^n$,
$n\geq 2$.  We present the results of the previous sections in a more
explicit form.  We obtain a sampling and interpolation theorem for the
weighted Bergman spaces on $\B$.  For each weight, the restriction to
$X$ induces an isomorphism from the weighted Bergman space on $\B$ to
the one on $X$ if $L$ is sufficiently ample.  This is in contrast to
Seip's result that no sequence in the disc is both sampling and
interpolating for any weighted Bergman space \cite{Sei}.  Also, every
bounded holomorphic function on $X$ extends to a unique function of just
barely exponential growth on $\B$, i.e., a function in the intersection
of all the nontrivial weighted Bergman spaces on $\B$, when $L$ is
sufficiently ample, for instance when $L$ is the $m$-th tensor power of
the canonical bundle $K$ with $m\geq 2$.  Whether the extension is
itself bounded is an important open question. 

In Section 5, assuming that the covering group is an arithmetic subgroup
of the automorphism group $PU(1,n)$ of $\B$, we establish two
dichotomies related to that in \cite{L\ac ar2} but using very different
means.  One of them says that either every holomorphic function $f$
continuous up to the boundary on the preimage of a $K^{\otimes m}$-curve
in a finite covering of $M$ extends to a continuous function on $\overline\B$
which is holomorphic on $\B$, or the boundary functions of such
functions $f$ generate a dense subspace of the space of continuous
functions $\d\B\to\C$.  This is probably another manifestation of the
elusive phenomenon discovered in \cite{L\ac ar2}.  Section 6 contains an
analogous dichotomy for harmonic functions.  In contrast to the case of
holomorphic functions, it is easy to see that a harmonic function on
$X$, continuous up to the boundary, generally does not extend to a
plurisubharmonic function bounded above on $\B$. 

\smallskip \flushpar {\it Acknowledgement.}  I would like to thank Bo
Berndtsson for helpful discussions.

\specialhead 1.  An extension theorem \endspecialhead

\flushpar
We will be working with the following objects.
\roster
\item  A covering space $\pi:Y\to M$ of a compact $n$-dimensional
K\"ahler manifold $M$ with a K\"ahler form $\omega$. 
\item A smooth function $\phi:Y\to\R$ such that $d\phi$ is bounded.
\item A line bundle $L$ on $M$ with canonical connection $\nabla$ and
curvature $\Theta$ in a hermitian metric $h$.
\item A section $s$ of $L$ over $M$ with $\nabla s\neq 0$ at each point of
its nonempty zero locus $C$.  Then $C$ is a smooth (possibly disconnected)
hypersurface in $M$.  Let $X=\pi^{-1}(C)$.
\endroster
We denote the pullbacks to $Y$ of $\omega$, $L$, and $s$ by the same
letters. 

Let $\O_\phi(X)$ be the vector space of holomorphic functions $f$ on
$X$ such that $f^2 e^{-\phi}$ is integrable on $X$ with respect to the
volume form of the induced K\"ahler metric on $X$.  This is a Hilbert
space with respect to the inner product 
$$(f,g)\mapsto\int_Xf\bar g e^{-\phi}\omega^{n-1}.$$  
We define $\O_\phi(Y)$ similarly. 

Let $U_1,\dots,U_m$ be the pullbacks of shrunk coordinate polydiscs
in which $C=\{z_n=0\}$, covering a neighbourhood of $C$. 
If $f\in\O_\phi(Y)$ and $x=(a_1,\dots,a_{n-1},0)
\in X\cap U_k$, then 
$$|f(x)|^2\leq c\int_D |f|^2,$$
where $D$ is the disc with $z_j=a_j$ for $j=1,\dots,n-1$ in $U_k$, and
$c>0$ is a constant independent of $f$.  Hence,
$$\int_{X\cap U_k}|f|^2e^{-\phi}\leq c\int_{U_k}|f|^2e^{-\phi},$$
so
$$\int_X|f|^2e^{-\phi}\leq c\int_Y|f|^2e^{-\phi}.$$
This shows that we have a continuous linear restriction map
$$\rho:\O_\phi(Y)\to\O_\phi(X), \qquad f\mapsto f|X.$$

In a previous paper we showed that under suitable curvature assumptions,
$\rho$ is surjective when $n\geq 2$.

\proclaim{1.1.  Theorem \cite{L\ac ar1, Thm\.  3.1}} If $n\geq 2$ and
$\Theta\geq i\d\bar\d\phi+\epsilon\omega$ for some $\epsilon>0$, then
$\rho$ is surjective.  
\endproclaim

By \cite{L\ac ar1, Cor\.  2.4}, since the weighted metric $e^\phi h$ in
$L$ has curvature $-i\d\bar\d\phi+\Theta\geq\epsilon\omega$, the $k$-th
$L^2$ cohomology group $H_{(2)}^k(Y,L^\vee)$ of $Y$ with coefficients in
the dual bundle $L^\vee$ with the dual metric $e^{-\phi}h^\vee$ vanishes
for $k<n$.  The proof of Theorem 1.1 is based on vanishing for $k=1$,
for which we need $n\geq 2$.  We will now use vanishing for $k=0$ to
show that $\rho$ is injective. 

Let $f\in\O_\phi(Y)$ such that $f|X=0$.  Then $\alpha=fs^\vee$ is a
holomorphic section of $L^\vee$ on $Y$.  We will show that $\alpha$ is
square-integrable with respect to $e^{-\phi}h^\vee$.  Then
vanishing of $H_{(2)}^0(Y,L^\vee)$ implies that $\alpha=0$, so $f=0$. 

Since $\nabla s\neq 0$ on $C$, there is a constant $c>0$ such that
$\text{dist}(\cdot,C)\leq c|s|$ on $M$.  For $y\in Y \setminus X$ let
$x\in X$ have $\text{dist}(y,x)=\text{dist}(y,X)$.  Then
$$\multline |\alpha(y)|=|f(y)||s(y)|^{-1}\leq c|f(x)-f(y)|/\text{dist}(x,y)
\\ \leq c\sup|df|\leq c\int_B |f|\leq c\Big(\int_B|f|^2\Big)^{1/2},
\endmultline $$
where the supremum is taken over a ball centred at $y$ covering all of $M$,
and $B$ is a ball of larger radius.  Hence,
$$|\alpha(y)|^2 e^{-\phi(y)} \leq c\int_B |f|^2 e^{-\phi},$$
so
$$\int_Y|\alpha|^2 e^{-\phi}\leq c\int_Y|f|^2e^{-\phi}<\infty.$$

We have proved the following theorem. 

\proclaim{1.2. Theorem}  Suppose 
$$\Theta\geq i\d\bar\d\phi +\epsilon\omega$$
for some $\epsilon>0$.  Then $\rho$ is injective.  If $\dim X\geq 1$, then
$\rho$ is an isomorphism.
\endproclaim

By induction, the theorem generalizes to the case when $C$ is the common
zero locus of sections $s_1,\dots,s_k$, $k\leq n$, of $L$ over $M$
which, in a trivialization, can be completed to a set of local
coordinates at each point of $C$.  When $k=n-1$, such $C$ will be
referred to as $L$-curves.  If $L$ is very ample, and therefore the
pullback of the hyperplane bundle by an embedding of $M$ into some
projective space, then this condition means that the linear subspace
$\{s_1,\dots,s_k=0\}$ intersects $M$ transversely in a smooth subvariety
$C$ of codimension $k$.  By Bertini's theorem, this holds for the
generic linear subspace of codimension $k$.  If $k\leq n-1$, then $C$ is
connected and the map $\pi_1(C)\to\pi_1(M)$ is surjective by the
Lefschetz hyperplane theorem, which implies that $X$ is connected. 

An important example of a function $\phi$ as above is obtained by
smoothing the distance $\delta$ from a fixed point in $Y$.  By a result
of Napier \cite{Nap}, there is a smooth function $\tau$ on $Y$ such that
\roster
\item $c_1\delta\leq\tau\leq c_2\delta+c_3$ for some $c_1, c_2, c_3>0$,
\item $d\tau$ is bounded, and
\item $i\d\bar\d\tau$ is bounded.
\endroster
Furthermore, by \therosteritem1 and since the curvature of $Y$ is
bounded below, there is $c>0$ such that $e^{-c\tau}$ is integrable on
$Y$.  Then $e^{-c\tau}$ is also integrable on $X$.

If $L$ is positive, then $k\Theta\geq i\d\bar\d\tau$ for $k\in\Bbb N$
sufficiently large by \therosteritem3, so the curvature inequality in
Theorem 1.2 holds if $L$ is replaced by a sufficiently high tensor power
of itself. 

\smallskip\flushpar
{\bf 1.3. Example.}  There is an open hyperbolic
Riemann surface $X$ such that for some $\epsilon>0$, any $f\in\O(X)$
with $|f|\leq ce^{\epsilon\delta}$ is constant, where $\delta$ is the
distance from a fixed point in the Poincar\ac e metric.

Namely, there is an example due to Cousin of a projective 2-dimensional
torus (abelian surface) $M$ with a $\Bbb Z$-covering space $Y\to M$ such
that $Y$ has no nonconstant holomorphic functions \cite{NR, 3.9}.  Let
$\tau$ be a smooth function on $Y$ satisfying \therosteritem1,
\therosteritem2, and \therosteritem3, such that $e^{-\tau/2}$ is
integrable.  Let $L$ be a very ample line bundle on $M$ such that
$\Theta\geq i\d\bar\d\tau +\epsilon\omega$.  Let $X$ be the pullback in
$Y$ of an $L$-curve in $M$.  Then $\O_\tau(X)=\C$ by Theorem 1.2.  If
$f\in\O(X)$ and $|f|\leq ce^{\epsilon\delta}$ with $\epsilon>0$
sufficiently small, then $|f|\leq ce^{\tau/4}$, so $f\in\O_\tau(X)$ and
$f$ is constant.

\specialhead 2.  New examples of Riemann surfaces with corona \endspecialhead

\flushpar Let $X$ be a complex manifold.  Let $H^\infty(X)$ be the space
of bounded holomorphic functions on $X$, which is a Banach algebra in
the supremum norm.  Let $\Cal M$ be the character space of
$H^\infty(X)$, which is a compact Hausdorff space in the weak-star
topology.  There is a continuous map $\iota:X\to\Cal M$ taking $x\in X$
to the evaluation character $f\mapsto f(x)$.  
The Corona Problem asks whether $\iota(X)$ is dense in $\Cal M$. 
The complement of the closure of $\iota(X)$ in $\Cal M$ is
referred to as the corona, so if $\iota(X)$ is not dense in $\Cal M$, 
then $X$ is said to have corona.  It is well known that the following 
are equivalent. 
\roster 
\item $\iota(X)$ is dense in $\Cal M$.  
\item If $f_1,\dots,f_m\in H^\infty(X)$ and $|f_1|+\dots+|f_m|>
\epsilon>0$, then there are $g_1,\dots,g_m\in H^\infty(X)$ such that 
$f_1 g_1+\dots+f_m g_m=1$.  
\endroster 
By Carleson's famous Corona Theorem (1962), the disc has no corona.  The
Corona Theorem holds for Riemann surfaces of finite type and planar
domains of various kinds.  The Corona Problem for arbitrary planar
domains is open.  Around 1970, Cole constructed the first example of a
Riemann surface with corona \cite{Gam}.  By modifying Cole's example,
Nakai obtained a regular Parreau-Widom surface with corona \cite{Nak},
\cite{Has, p\. 229}.  Recently, Barrett and Diller showed that the
homology covering spaces of domains in the Riemann sphere, whose
complement has positive logarithmic capacity and zero length, have
corona \cite{BD}.  See also \cite{EP, 7.3}.  Sibony \cite{Sib} found the
first example of a domain of holomorphy in $\C^n$, $n\geq 2$, with
corona.  There are no known examples of such domains without corona. 

We will now present a new class of Riemann surfaces with corona (see
also Theorem 4.2).  We remind the reader that if $Y$ is a bounded domain
in $\C^n$ covering a compact complex manifold $M$, then $Y$ is a domain
of holomorphy \cite{Sie, p\.  136} and $M$ is projective.  In fact, the
canonical bundle of $M$ is ample \cite{Kol, 5.22}. 

\proclaim{2.1.  Theorem} Let $\pi:Y\to M$ be a covering map, where $Y$
is a bounded domain in $\C^n$, $n\geq 2$, and $M$ is compact.  Let $L$
be an ample line bundle on $M$.  If $C$ is an $L^{\otimes m}$-curve in
$M$ with $m$ sufficiently large, then the Riemann surface
$X=\pi^{-1}(C)$ has corona.  In fact, the natural map from $X$ into the
character space $\Cal M$ of $H^\infty(X)$ extends to an embedding of
$Y$ into $\Cal M$ which maps $Y\setminus X$ into the corona of $X$. 

\endproclaim

\demo{Proof} Let $\tau\geq 0$ be a smoothed distance function on $Y$ as
described in Section 1, such that $e^{-\tau}$ is integrable on $Y$, and
hence on $X$, so $H^\infty(X)\subset\O_\tau(X)$.   

We claim that $\O_\tau(X)\cdot\O_\tau(X)\subset\O_{3\tau}(X)$.  Namely,
suppose $f\in\O_\tau(X)$.  For $p\in X$, let $B$ be the ball of
radius 1 centred at $p$ in a metric pulled up from $C$.  Then 
$$|f(p)|\leq c\int_B |f|\leq c\big(\int_B |f|^2 \big)^{1/2},$$
where the constants are independent of $p$, so
$$|f(p)|e^{-\tau(p)/2} \leq c\big(\int_B|f|^2 e^{-\tau}\big)^{1/2}
<\infty.$$
Hence, $|f|\leq ce^{\tau/2}$, so 
$$\int_X |f|^4 e^{-3\tau} \leq c \int_X e^{2\tau} e^{-3\tau}<\infty,$$
and $f^2\in\O_{3\tau}(X)$.

If $m$ is sufficiently large, then the curvature of $L^{\otimes m}$
is at least $3i\d\bar\d\tau+\omega$
for some K\"ahler form $\omega$ on $M$, so the restriction map
$\rho:\O_{3\tau}(Y)\to \O_{3\tau}(X)$ is an isomorphism by Theorem 1.2. 
Clearly, $\rho\O_\tau(Y)=\O_\tau(X)$.
For $p\in Y$, let $\lambda_p$ be the linear functional
$f\mapsto\rho^{-1}(f)(p)$ on $H^\infty(X)$.  If $f,g\in H^\infty(X)$,
then $\rho^{-1}(f)\rho^{-1}(g)\in\O_\tau(X)\cdot\O_\tau(X)\subset
\O_{3\tau}(X)$, and $\rho^{-1}(f)\rho^{-1}(g)|X=fg$, so
$\rho^{-1}(f)\rho^{-1}(g)=\rho^{-1}(fg)$.  This shows that $\lambda_p$
is a character on $H^\infty(X)$. 

We have obtained a map $\iota:Y\to\Cal M$, $p\mapsto\lambda_p$, extending
the natural map from $X$ into $\Cal M$.  We claim that $\iota$ is a
homeomorphism onto its image with the induced topology, and that
$\overline{\iota(X)}\cap\iota(Y)=\iota(X)$.  First of all, since
$H^\infty(Y)$ separates points, $\iota$ is injective.  The topology on
$\Cal M$ is the weakest topology that makes all the Gelfand transforms
$\hat f:\chi\mapsto\chi(f)$, $f\in H^\infty(X)$, continuous, so $\iota$
is continuous if and only if $\hat f\circ\iota$ is continuous for all
$f\in H^\infty(X)$, but $\hat f\circ\iota=\rho^{-1}(f)$, so this is
clear.  Let $p\in Y$ and suppose the polydisc $P=\{z\in\C^n:
|f_i(z)|<\epsilon, i=1,\dots,n\}$, where $f_i(z)=z_i-p_i$, is contained
in $Y$.  Let $V=\{\chi\in\Cal M:|\chi(f_i)|<\epsilon, i=1,\dots,n\}$. 
Then $V$ is open in $\Cal M$, and $\iota(P)=V\cap\iota(Y)$.  This shows
that $\iota:Y\to\iota(Y)$ is open, so $\iota$ is a homeomorphism onto
its image.  Finally, suppose $p\in Y\setminus X$.  Say
$\max|f_i|>\epsilon>0$ on $X$.  Then $V$ is an open neighbourhood of
$\iota(p)$ which does not intersect $\iota(X)$.  Hence,
$\overline{\iota(X)}\cap\iota(Y)=\iota(X)$.  This shows that the corona
of $X$ contains an embedded image of $Y\setminus X$.  \qed\enddemo

Let us remark that a Riemann surface as in the theorem is not
Parreau-Widom, since a Parreau-Widom surface $X$ embeds into the
character space of $H^\infty(X)$ as an open subset \cite{Has, p\. 222}. 
See also the proof of Theorem 5.2.

By the same argument we easily obtain the following more general result.

\proclaim{2.2.  Theorem} Let $Y$ be a covering space of a projective 
manifold $M$ with $\dim M\geq 2$, $L$ be a line bundle on $M$, and $C$ 
be an $L$-curve in $M$ with preimage $X$ in $Y$.  If
\roster 
\item  $L$ is sufficiently ample, and
\item  there is a bounded holomorphic map $g:Y\to\C^m$ with
$\overline{g(X)}\neq g(Y)$,
\endroster
then $X$ is a Riemann surface with corona.  
\endproclaim

Finally, under the hypotheses of Theorem 2.1, if $U\subset Y$ is a
neighbourhood of $X$ with $Y\not\subset\overline U$, then we can show by
the same argument as in the proof of Theorem 2.1 that $U$ has corona. 
Since $X$ is Stein, $U$ may be chosen to be a domain of holomorphy by
Siu's theorem \cite{Siu}.  Thus we obtain new examples of bounded
domains of holomorphy in $\C^n$ with corona.

\specialhead 3.  Generating H\"ormander algebras on covering spaces
\endspecialhead

In this section, we adapt results of H\"ormander \cite{H\"or} on
generating algebras of holomorphic functions of exponential growth to
the case of covering spaces over compact manifolds.  We let $X\to M$
be a covering space of an $n$-dimensional compact hermitian manifold $M$.
Let $\phi:X\to[0,\infty)$ be a smooth function such that
\roster
\item $d\phi$ is bounded, and
\item $e^{-c\phi}$ is integrable on $X$ for some $c>0$.
\endroster
Let
$$\A_\phi=\A_\phi(X)=\bigcup_{c>0}\O_{c\phi}(X)$$
be the vector space of holomorphic functions $f$ on $X$ such that 
$f^2e^{-c\phi}$ is integrable on $X$ for some $c>0$.  By
\therosteritem2, $\A_\phi$ contains all bounded holomorphic functions
on $X$.  The following is easy to see by an argument similar to that in
the proof of Theorem 2.1.

\proclaim{3.1. Proposition}  A holomorphic function $f$ on $X$ is in
$\A_\phi$ if and only if $|f|\leq ce^{a\phi}$ for some $a>0$.
Hence, $\A_\phi$ is a $\C$-algebra, called a H\"ormander algebra.  
\endproclaim

If functions $f_1,\dots,f_m$ in $\A_\phi$ generate $\A_\phi$, then there
are $g_1,\dots,g_m\in\A_\phi$ such that $f_1 g_1+\dots+f_m g_m=1$, so
$$\max_{i=1,\dots,m}|f_i|\geq ce^{-a\phi} \qquad\text{for some }a>0.$$
We will establish an effective converse to this observation.  Our proof
is a straightforward adaptation of H\"ormander's Koszul complex argument
in \cite{H\"or}.  See also \cite{EP, 7.3}.

Let $m\geq 1$ and $r,s\geq 0$ be integers, and $t\in[0,\infty)$.  Choose
a basis $\{e_1,\dots,e_m\}$ for $\C^m$.  Let $L_r^s(t)$ be the space of
smooth $\Lambda^s\C^m$-valued $(0,r)$-forms on $X$ which are
square-integrable with respect to $e^{-t\phi}$.  

\proclaim{3.2.  Lemma} Suppose $M$ is K\"ahler with K\"ahler form
$\omega$, and $\Ric(X)+ti\d\bar\d\phi\geq \epsilon\omega$ for some
$\epsilon>0$.  If $\eta\in L_{r+1}^s(t)$ and $\bar\d\eta=0$, then there
is $\xi\in L_r^s(t)$ with $\bar\d\xi=\eta$.  
\endproclaim

\demo{Proof}  This follows directly from standard $L^2$ theory.  See for
instance \cite{Dem, Sec\. 14}.
\qed\enddemo

Now let $f_1,\dots,f_m$ be holomorphic functions on $X$ such that
$$ce^{-c_2\phi}\leq\max_i|f_i|^2\leq ce^{c_1\phi}, \qquad c_1,c_2>0.$$
Define a linear operator $\alpha:L_r^{s+1}(t)\to L_r^s(t+c_1)$ by the
formula 
$$\alpha(e_{i_1}\wedge\dots\wedge e_{i_{s+1}})=\sum_{k=1}^{s+1}(-1)^{k+1}
f_{i_k} e_{i_1}\wedge\dots\wedge\hat e_{i_k}\wedge\dots\wedge e_{i_{s+1}},$$
and set $\alpha=0$ on $L_r^0(t)$.  Then $\alpha^2=0$, and $\alpha$
commutes with $\bar\d$.  Also define a linear operator
$\beta:L_r^s(t)\to L_r^{s+1}(t+c_1+2c_2)$ by the formula
$$\beta(\xi)=\frac1{|f_1|^2+\dots+|f_m|^2}\big(\sum_{i=1}^m\bar f_i e_i
\big)\wedge\xi.$$
Then $\alpha\beta+\beta\alpha$ is the inclusion $L_r^s(t)\hookrightarrow
L_r^s(t+2(c_1+c_2))$. 

\proclaim{3.3. Lemma}  Suppose $e^{-a\phi}$ is integrable on $X$.  If
$\xi\in L_r^s(t)$ and $\alpha(\xi)=0$, then there is $\eta\in
L_r^{s+1}(t+c_1+2c_2)$ such that $\alpha(\eta)=\xi$ and in addition
$\bar\d\eta\in L_{r+1}^{s+1}(t+2c_1+3c_2+a)$ if $\bar\d\xi=0$.
\endproclaim

\demo{Proof}  Take $\eta=\beta(\xi)$.  Say $\bar\d\xi=0$.  Then
$$|\bar\d\eta|\leq c|\xi|\max_i \big|\bar\d\frac{\bar f_i}{|f_1|^2+\dots+
|f_m|^2}\big|\leq c|\xi|\max_i|\d f_i|e^{(c_1+3c_2)\phi/2}.$$
If $f\in\O(X)$ and $|f|^2\leq ce^{c_1\phi}$, then $|\d f|^2\leq c
e^{(c_1+a)\phi}$.  Hence, $|\bar\d\eta|^2\leq c|\xi|^2e^{(2c_1+3c_2+a)\phi}$,
and the lemma follows.
\qed\enddemo

\proclaim{3.4. Theorem}  Let $X\to M$ be a covering space of a
compact K\"ahler manifold with K\"ahler form $\omega$.  Let 
$\phi:X\to[0,\infty)$ be a smooth plurisubharmonic function such that
\roster
\item $d\phi$ is bounded,
\item $e^{-a\phi}$ is integrable on $X$, $a>0$, and
\item $\Ric(X)+bi\d\bar\d\phi\geq\epsilon\omega$ for some $\epsilon>0$.
\endroster
Let $f_1,\dots,f_m$ be holomorphic functions on $X$ such that
$$ce^{-c_2\phi}\leq\max_i|f_i|^2\leq ce^{c_1\phi}, \qquad c_1,c_2>0.$$
Let $r,s\geq 0$ be integers.  If $\xi\in L_r^s(t)$, 
$\bar\d\xi=0$, $\alpha(\xi)=0$, and $t\geq b-2c_2$, then there is 
$\eta\in L_r^{s+1}(u)$, where
$$u=t+c_1+2c_2+(m-s-1)(a+3(c_1+c_2)),$$
such that $\bar\d\eta=0$ and $\alpha(\eta)=\xi$.
\endproclaim

Taking $r=s=0$ and $\xi=1$, we obtain the following corollary.

\proclaim{3.5. Corollary}  Let the hypotheses be as in Theorem 3.4.  
There are
$g_1,\dots,g_m$ in $\O_{t\phi}(X)$ where 
$$t=\max\{a,b-2c_2\}+c_1+2c_2+(m-1)(a+3(c_1+c_2)),$$
such that $f_1 g_1+\dots+f_m g_m=1$.
\endproclaim

The following corollary is our analogue of \cite{H\"or, Thm\. 1}.

\proclaim{3.6. Corollary}  Let $X\to M$ be a covering space of a
compact K\"ahler manifold with K\"ahler form $\omega$.  Let
$\phi:X\to[0,\infty)$ be a smooth function such that
\roster
\item $d\phi$ is bounded,
\item $e^{-c\phi}$ is integrable on $X$ for some $c>0$, and
\item $i\d\bar\d\phi\geq\epsilon\omega$ for some $\epsilon>0$.
\endroster
Then functions $f_1,\dots,f_m$ in $\A_\phi$ generate $\A_\phi$ if and only
if 
$$\max_{i=1,\dots,m}|f_i|\geq ce^{-a\phi} \qquad\text{for some }a>0.$$
\endproclaim

The hypotheses of the corollary are satisfied for example when $X$ is
the unit ball in $\C^n$ and $\phi=-\log(1-|\cdot|^2)$, which is
comparable to the Bergman distance from the origin (see Section 4).

\demo{Proof of Theorem 3.4} If $s\geq m$ or $r>n$, then $\xi=0$ and we
take $\eta=0$.  Assume that $s<m$ and $r\leq n$, and that the theorem
has been proved with $r$, $s$ replaced by $r+1$, $s+1$.  By Lemma 3.3,
there is $\eta_1\in L_r^{s+1}(t+c_1+2c_2)$ such that $\alpha(\eta_1)=\xi$
and $\bar\d\eta_1\in L_{r+1}^{s+1}(t+2c_1+3c_2+a)$.  Now
$\bar\d\bar\d\eta_1=0$ and $\alpha(\bar\d\eta_1)= \bar\d\alpha(\eta_1)
=\bar\d\xi=0$, so by the induction hypothesis, there is $\eta_2\in
L_{r+1}^{s+2}(u-c_1)$ such that $\bar\d\eta_2=0$ and
$\alpha(\eta_2)=\bar\d\eta_1$.  By Lemma 3.2, there is $\eta_3\in
L_r^{s+2}(u-c_1)$ such that $\bar\d\eta_3=\eta_2$.  Now let
$\eta=\eta_1-\alpha(\eta_3)\in L_r^{s+1}(u)$.  Then
$\bar\d\eta=\bar\d\eta_1-\alpha(\bar\d\eta_3)
=\bar\d\eta_1-\alpha(\eta_2)=0$ and $\alpha(\eta)=\alpha(\eta_1)=\xi$. 
\qed\enddemo

\specialhead 4.  The case of the ball  \endspecialhead

In this section, we will consider the results of the previous sections
in the explicit setting of the unit ball $\B$ in $\C^n$, $n\geq 2$.  For
an instructive discussion of compact ball quotients, see \cite{Kol, Ch\. 
8}. 

Let $M$ be a projective manifold covered by $\B$ with a positive line
bundle $L$ with curvature $\Theta$.  Let $X$ be the preimage in $\B$ of
an $L$-curve $C$ in $M$.  We are particularly interested in the
extension problem for bounded holomorphic functions on $X$. 

The restriction map $H^\infty(\B)\to H^\infty(X)$ is injective, so we
can consider $H^\infty(\B)$ as a subspace of $H^\infty(X)$, which is
closed in the sup-norm.  In the locally uniform topology, however,
$H^\infty(\B)$ is dense in $H^\infty(X)$.  Namely, say $f\in
H^\infty(X)$, and let $F\in\O(\B)$ be an extension of $f$.  For $r<1$,
the functions $z\mapsto F(rz)$ are bounded in $\B$, and they converge
locally uniformly to $f$ on $X$ as $r\to 1$.

It is well known that if $Y$ is a complex submanifold of a neighbourhood
of $\overline\B$, then every bounded holomorphic function on $Y\cap\B$
extends to a bounded holomorphic function on $\B$ \cite{HL, 4.11.1}. 
Our $X$ is of course far from extending to a submanifold of a larger
ball, and the extension problem for bounded holomorphic functions on $X$
is very much open.  We will show, however, that when $L$ is sufficiently
positive, bounded holomorphic functions on $X$ extend to holomorphic
functions on $\B$ whose growth is, in a precise sense, just barely
exponential. 

First we collect a few formulas concerning the geometry of $\B$.  The
Bergman
metric of $\B$ is 
$$ds^2=\sum_{j,k=1}^n g_{jk}dz_j\otimes d\bar z_k,$$
where
$$g_{jk}=\langle \frac{\d}{\d z_j},\frac{\d}{\d z_k}\rangle =
\frac{n+1}{(1-|z|^2)^2} ((1-|z|^2)\delta_{jk} + \bar z_j z_k).$$
This is a common convention.  It is used for instance in \cite{Kra} and
\cite{Sto}; other authors may use a constant scalar multiple of the
above. 
The K\"ahler form of the Bergman metric is
$$\omega=-\frac 1 2 \Im ds^2=\frac i 2 \sum_{j,k=1}^n g_{jk} dz_j
\wedge d\bar z_k=-\frac {i(n+1)} 2\d\bar\d\log(1-|z|^2).$$
The distance from the origin in the Bergman metric is
$$\delta(z)=\frac{\sqrt{n+1}}2\log\frac{1+|z|}{1-|z|},\qquad z\in\B.$$
The Ricci curvature of the Bergman metric, i.e., the curvature form of
the induced metric in the cocanonical bundle (the top exterior power of
the
tangent bundle) is
$$\Ric(\omega)=-\frac i 2 \d\bar\d\log\det(g_{jk})=-\omega.$$

Let 
$$\tau(z)=-\frac{n+1} 2\log(1-|z|^2).$$  
Then $i\d\bar\d\tau=\omega$, and $\tau$ is a nonnegative strictly 
plurisubharmonic exhaustion of $\B$. 
Also, $\tau$ is comparable to $\delta$.  More precisely,
$$\sqrt{n+1}\,\delta\leq\tau\leq\sqrt{n+1}\,\delta+(n+1)\log 2.$$
It may be shown that $|d\tau|$ is bounded.  
Let $\Omega=\omega^n/n!$ be the volume form of the Bergman metric.  Then
$$\Omega=c(1-|z|^2)^{-(n+1)}\Omega_0,$$ where $c>0$ is constant, and
$\Omega_0$ is the euclidean volume form on $\C^n$.  We have
$$\int_\B e^{-c\tau}\Omega<\infty\quad\text{if and only if}\quad
c>\dfrac{2n}{n+1}.$$  

The weighted Bergman space $\Cal B_p^a$ is the space of holomorphic
functions $g$ on $\B$ such that  
$$\int_\B|g|^p (1-|z|)^a\omega^n<\infty.$$
See [Sto, Ch\. 10].  We have $\Cal B_2^n=0$.  The
intersection $\Cal B_2^{n+0}=\bigcap\limits_{\epsilon>0}\Cal
B_2^{n+\epsilon}$ contains the Hardy space $H^2(\B)$ of
holomorphic functions $f$ on $\B$ such that $|f|^2$ has a harmonic
majorant. 

The boundary behaviour of the functions in $\Cal B_2^{n+0}$ may be
rather wild.  A theorem of Bagemihl, Erd\"os, and Seidel \cite{BES},
\cite{Mac}, states that if $\mu:[0,1)\to[0,\infty)$ goes to infinity at
$1$, then there exists a holomorphic function $f$ on the unit disc with
$|f|\leq \mu(|\cdot|)$, such that for some sequence $r_n\nearrow 1$, we
have $\min\limits_{|z|=r_n} |f(z)|\to\infty$.  In particular, $f$ does
not have a finite limit along any curve that intersects every
neighbourhood of the boundary.  Taking $\mu(r)=-\log(1-r)$, we get $f$
in $\Cal B_2^{1+0}$. 
  
We have
$$\Cal B_2^a=\O_{\frac{2a}{n+1}\tau}(\B).$$
Let
$$\Cal E(\B)=\bigcap_{c>\frac{2n}{n+1}}\O_{c\tau}(\B),$$
and define $\Cal E(X)$ similarly.  These are Fr\ac echet-Hilbert spaces.

From now on we let $n\geq 2$.  Theorem 1.2 yields the following result. 

\proclaim{4.1. Theorem}  If $\Theta>c\omega$, then the restriction map
$\rho:\O_{c\tau}(\B)\to\O_{c\tau}(X)$ is an isomorphism.

Suppose 
$$\Theta>\frac{2n}{n+1}\,\omega.$$
This holds for instance if $L=K^{\otimes k}$ with $k\geq 2$.  Then we
have a well-defined restriction map $\rho:\Cal E(\B)\to\Cal E(X)$, which
is an isomorphism.  Hence, every bounded holomorphic function on $X$
extends to a unique function in $\Cal E(\B)$.
\endproclaim

The 1-dimensional case is very different.  Seip has shown that no
submanifold of the disc is both sampling and interpolating for any
weighted Bergman space \cite{Sei}. 

Now we can easily obtain fairly explicit examples of Riemann surfaces
with corona.

\proclaim{4.2. Theorem}  If $\Theta>\dfrac{2n}{n+1}\,\omega$, then $X$ has
corona.
\endproclaim

\demo{Proof}  Let $p\in\B\setminus X$ and $f_i(z)=z_i-p_i$,
$i=1,\dots,n$.  Then $\sum|f_i|>\epsilon>0$ on $X$.  Suppose $X$
has no corona.  Then there are $g_1,\dots,g_n\in H^\infty(X)$ with
$\sum f_i g_i=1$.  By Theorem 4.1, $g_i$ extends to a function
$G_i\in\Cal E(\B)$.  Then $h=\sum f_i G_i\in\Cal E(\B)$, and $h|X=1$. 
Again by Theorem 4.1, $h=1$, which is absurd since $h(p)=0$.
\qed\enddemo

The proof of Theorem 2.1 shows that if $\Theta>\dfrac{6n}{n+1}\omega$,
then $\B$ embeds into the character space $\Cal M$ of $H^\infty(X)$,
taking $\B\setminus X$ into the corona.  

\smallskip
Now $\tau$ satisfies the hypotheses of Corollary 3.6, both on $\B$ and
on $X$.  Hence, functions $f_1.\dots,f_m\in\A_\tau(\B)$ generate
$\A_\tau(\B)$ if and only if 
$$\max_{i=1,\dots,m}|f_i|\geq ce^{-a\tau} \qquad\text{for some }a>0,$$ 
and similarly for $X$. 

\smallskip
In contrast to the first part of Theorem 4.1, we we obtain the following
result from Corollary 3.5. 

\proclaim{4.3.  Theorem}  If $s\geq 1$ and $\Theta\leq s\omega$, then the
restriction map $\rho:\O_{c\tau}(\B)\to\O_{c\tau}(X)$ is not an
isomorphism for 
$$c>s+\frac{2n^2-n+1}{n+1}.$$
\endproclaim

\demo{Proof} By adjunction, the Ricci curvature of $X$ in the metric
induced by the Bergman metric on $\B$ is $\Ric(X)=-\Theta-\omega$.  We
apply Corollary 3.5 with $X\to C$, $\phi=\tau$,
$a=2n/(n+1)+\epsilon$, $b=s+1+\epsilon$, $\epsilon>0$, $f_i=z_i-p_i$,
$i=1,\dots,n$, where $p\in\B\setminus X$, and $c_1,c_2=0$.  For every
$c>s+(2n^2-n+1)/(n+1)$ we obtain $g_1,\dots,g_n$ in $\O_{c\tau}(X)$ such
that $f_1 g_1+\dots+f_n g_n=1$.  If
$\rho:\O_{c\tau}(\B)\to\O_{c\tau}(X)$ was an isomorphism, there would be
$G_1,\dots,G_n$ in $\O_{c\tau}(\B)$ with $\sum(z_i-p_i)G_i=1$ on $\B$,
which is absurd.  
\qed\enddemo

Consider the special case when $L=K^{\otimes m}$, $m\geq 1$.  We have
proved that $\rho:\O_{c\tau}(\B)\to\O_{c\tau}(X)$ is an isomorphism if
$c<m$, but not if $c>m+\dfrac{2n^2-n+1}{n+1}$.  Furthermore, we can
easily show that $\rho$ is not injective if $c>2m+\dfrac{2n}{n+1}$. 

Namely, $e=dz_1\wedge\dots\wedge dz_n$ is a zero-free section of $K$
with norm $ce^{-\tau}$.  Let $s\neq 0$ be a holomorphic section of $L$
on $M$, vanishing on $C$.  Then $f=s/e^{\otimes m}$ is a holomorphic
function in the kernel of $\rho$ with $|f|\leq ce^{m\tau}$, so
$f\in\O_{c\tau}(\B)$ for all $c>2m+\dfrac{2n}{n+1}$.

\specialhead 5.  A dichotomy  \endspecialhead

\flushpar As before, we consider a projective manifold $M$ covered by
the unit ball $\B$ in $\C^n$, $n\geq 2$, with a positive line bundle
$L$, and the preimage $X$ in $\B$ of an $L$-curve $C$ in $M$.  We will
denote the covering group by $\G$.  The bounded extension problem for
holomorphic functions is related to the question of which bounded
harmonic functions on $X$ are real parts of holomorphic functions.  This
question was studied in \cite{L\ac ar2}, where the following dichotomy
was established in the more general setting of a nonelementary
Gromov-hyperbolic covering space of a compact K\"ahler manifold. 

\proclaim{5.1. Theorem \cite{L\ac ar2, Thm\. 4.2}}  One of the following
holds.
\roster
\item Every positive harmonic function on $X$ is the real part of a
holomorphic function.
\item If $u\geq 0$ is the real part of an $H^1$ function on $X$, then
the boundary decay of $u$ at a zero on the Martin boundary of $X$ is no
faster than its radial decay.
\endroster
\endproclaim

By results of Ancona \cite{Anc}, the Martin compactification of $X$ is
naturally homeomorphic to $X\cup\S$, where $\S=\d\B$ is the unit sphere.

Clearly, if \therosteritem1 holds, then there are holomorphic functions
on $X$ with a bounded real part that do not extend to a holomorphic
function on $\B$ with a bounded real part.  

If \therosteritem1 holds, then each Martin function $k_p$, $p\in\S$, is
the real part of a holomorphic function $f_p$ on $X$.  Then the
holomorphic map $\exp(-f_p):X\to\D$ is proper at every boundary point
except $p$.  Here, $\D$ denotes the unit disc.  Also, if $p,q\in\S$,
$p\neq q$, then the holomorphic map
$(\exp(-f_p),\exp(-f_q)):X\to\D\times\D$ is proper.  However, we have
the following result. 

\proclaim{5.2. Theorem}  There is no proper holomorphic map $X\to\D$.
\endproclaim

\demo{Proof} Bounded holomorphic functions separate points on $X$, so if
there is a proper holomorphic map $X\to\D$, then $X$ is Parreau-Widom by
a theorem of Hasumi \cite{Has, p\.  209}.  By \cite{L\ac ar2, Thm\. 
5.1}, if $X$ is Parreau-Widom, then $X$ is either isomorphic to $\D$ or
homeomorphic to the 2-sphere with a Cantor set removed.  Both
possibilities are excluded by the Martin boundary of $X$ being $\S$. 
\qed\enddemo

When $L$ is sufficiently ample, we can prove a stronger result.

\proclaim{5.3. Theorem}  If $L$ is sufficiently ample and $f$ is a
holomorphic function on $X$, then $f^{-1}(U)$ is not relatively compact
in $X$ for any nonempty open subset $U$ of the image $f(X)$.  In other
words, every value of $f$ is taken at infinity.
\endproclaim

\demo{Proof} Suppose there is a holomorphic function $f$ on $X$ such
that $f^{-1}(U)$ is relatively compact in $X$ for some nonempty open
subset $U$ of $f(X)$.  We may assume that $0\in U$.  Then $1/f$ is a
meromorphic function on $X$ which has a pole $p$ and is bounded outside
the compact closure of $f^{-1}(U)$.  Since bounded holomorphic functions
separate points on $X$, a theorem of Hayashi \cite{Hay} now implies that
the natural map from $X$ into the character space of $H^\infty(X)$ is
open when restricted to some neighbourhood of $p$.  By Theorem 2.1, this
is absurd when $L$ is sufficiently ample.
\qed\enddemo

We will now present another dichotomy in a similar vein.  Let $\Cal
C_K(\S)$ denote the space of continuous functions $\S\to K$, $K=\R$ or
$K=\C$, with the supremum norm.  Let $P$ be the subspace of $\Cal
C_\R(\S)$ of boundary values of pluriharmonic functions on $\B$ which
are continuous on $\overline\B$, and $\Cal O$ be the subspace of $\Cal
C_\C(\S)$ of boundary values of holomorphic functions on $\B$ which are
continuous on $\overline\B$.  It is known that if $V$ is a proper closed
subspace of $\Cal C_\R(\S)$ and $V$ is invariant under the action of the
automorphism group $G=PU(1,n)$ of $\B$, then $V=\R$ or $V=P$.  Also, if
$V$ is a proper closed $G$-invariant subspace of $\Cal C_\C(\S)$, then
$V$ is one of the following: $\C$, $\O$, $\overline\O$, $P+iP$
\cite{Rud, 13.1.4}. 

If $C$ is an $L$-curve in a finite covering of $M$ with preimage $X$ in
$\B$, then we denote by $E(C)$ the space of functions $\alpha\in\Cal
C_\R(\S)$ such that the harmonic extension $H[\alpha]=H_X[\alpha]$ of
$\alpha$ to $X$ is the real part of a holomorphic function on $X$. 
Clearly, $P\subset E(C)$.  Now $\alpha\in E(C)$ if and only if all the
periods of $H[\alpha]$ vanish, so $E(C)$ is closed in $\Cal C_\R(\S)$. 

\proclaim{5.4. Theorem}  Suppose that the covering group $\G$ is 
arithmetic, and that $L$
is a tensor power of the canonical bundle.  Then one (and only one)
of the following holds.
\roster
\item  $E(C)=P$ for every $L$-curve in a finite covering of $M$.
\item  The subspace of $\Cal C_\R(\S)$ generated by $E(C)$ for all
$L$-curves $C$ in finite coverings of $M$ is dense in $\Cal C_\R(\S)$.
\endroster
\endproclaim

Note that \therosteritem2 holds if \therosteritem1 in Theorem 5.1 holds
for the preimage $X$ of some $L$-curve in a finite covering of $M$.

\demo{Proof} Let $C$ be an $L$-curve in a finite covering $M_1$ of $M$,
with preimage $X$ in $\B$.  Then $M_1=\B/\G_1$, where $\G_1$ is a
subgroup of finite index in $\G$.  Let $g$ be an element of the
commensurability subgroup $\Comm(\G)$ in $G$.  This means that $\G$ and
$g\G g^{-1}$ are commensurable, i.e., their intersection is of finite
index in both of them.  Then $\G_2=\G_1\cap g\G_1 g^{-1}$ is a subgroup
of finite index in $\G_1$.  If $\alpha\in E(C)$, so $H_X[\alpha]=\Re f$
with $f$ holomorphic on $X$, then $f\circ g$ is holomorphic on $g^{-1}X$
and $H_{g^{-1}X}[\alpha\circ g]=\Re f \circ g$.  If $\g\in\G_2$, then
$\g=g^{-1}\g_1 g$ for some $\g_1\in\G_1$, so $\g g^{-1}X=g^{-1}\g_1
gg^{-1}X=g^{-1}\g_1 X=g^{-1}X$. 

Hence, $g^{-1}X$ is $\G_2$-invariant, so $g^{-1}X$ is the preimage of an
$L$-curve $C'$ in the finite covering $\B/\G_2$ of $M$ (here is where we
use the assumption that $L$ is a tensor power of the canonical bundle),
and $\alpha\circ g\in E(C')$. 

This shows that the subspace $E$ of $\Cal C_\R(\S)$ described in
\therosteritem2 is invariant under $\Comm(\G)$.  Since $\G$ is
arithmetic, $\Comm(\G)$ is Hausdorff-dense in $G$ \cite{Zim, 6.2.4} (and
in fact conversely), so the closure $\overline E$ of $E$ is a $G$-invariant
subspace of $\Cal C_\R(\S)$.  Hence, $\overline E$ is either $P$ or $\Cal
C_\R(\S)$, and the theorem follows.  
\qed\enddemo

If the spaces $E(C)$ are rigid in the sense that they do not change when
$C$ is varied in its linear equivalency class, then the theorem yields a
strong dichotomy.

\proclaim{5.5. Corollary}  Suppose that $\G$ is arithmetic, and that $L$
is a tensor power of the canonical bundle.  Suppose also that if $C_1$
and $C_2$ are $L$-curves in the same finite covering of $M$, then
$E(C_1)=E(C_2)$.  Then $E(C)$ is either $P$ or $\Cal C_\R(\S)$ for every
$L$-curve $C$ in a finite covering of $M$.
\endproclaim

We obtain analogous results for holomorphic functions.  If $C$ is an
$L$-curve in a finite covering of $M$ with preimage $X$ in $\B$, let us
denote by $F(C)$ the closed subspace of functions $\alpha\in\Cal
C_\C(\S)$ that extend to a holomorphic function on $X$.  Clearly,
$\O\subset F(C)$, but $F(C)$ is considerably smaller than $\Cal
C_\C(\S)$. 

\proclaim{5.6. Lemma}  $F(C)\cap (P+i\,\Cal C_\R(\S))=\O.$
\endproclaim

\demo{Proof} Let $\alpha\in F(C)\cap(P+i\,\Cal C_\R(\S))$, so
$H[\alpha]=f\in\O(X)$ and there is $u\in\Cal C_\R(\overline\B)$ such
that $u|\B$ is pluriharmonic and $u|\S=\Re\alpha$, so $u|X=\Re f$. 
There is $F\in\O(\B)$ such that $u=\Re F$ and $F|X=f$.  We need to show
that $F$ extends continuously to $\overline\B$. 

Now $F$ maps $\B$ into a vertical strip.  Let $\sigma$ be an isomorphism
from a neighbourhood of the closure of this strip in the Riemann sphere
onto $\D$.  Then $\sigma\circ F$ is a bounded holomorphic function on
$\B$ and $\sigma\circ F|X=\sigma\circ f$, so $\sigma\circ F$ has the
same nontangential boundary function $\sigma\circ\alpha$ as $\sigma\circ
f$.  Since $\sigma\circ\alpha$ is continuous, $\sigma\circ F$ extends
continuously to $\overline\B$, and so does $F=\sigma^{-1}\circ\sigma\circ F$.  
\qed\enddemo

\proclaim{5.7. Theorem}  Suppose that $\G$ is arithmetic, and that $L$
is a tensor power of the canonical bundle.  Then one of the 
following holds.
\roster
\item  $F(C)=\O$ for every $L$-curve in a finite covering of $M$.
\item  The subspace of $\Cal C_\C(\S)$ generated by $F(C)$ for all
$L$-curves $C$ in finite coverings of $M$ is dense in $\Cal C_\C(\S)$.
\endroster
Suppose furthermore that if $C_1$ and $C_2$ are $L$-curves in the same 
finite covering of $M$, then $F(C_1)=F(C_2)$.  Then $F(C)=\O$ for every
$L$-curve $C$ in a finite covering of $M$.
\endproclaim

\demo{Proof} By the same argument as in the proof of Theorem 5.4, the
closure of the subspace $F$ of $\Cal C_\C(\S)$ described in
\therosteritem2 is $G$-invariant, so it is $\O$, $P+iP$, or $\Cal
C_\C(\S)$ itself.  Lemma 5.6 shows that if $F\subset P+iP$, then $F=\O$. 
\qed\enddemo

Loosely speaking, either the spaces $F(C)$ are all the same, and equal
to $\O$, for all $L$-curves in finite coverings of $M$, or they are
diverse enough that every complex continuous function on $\S$ can be
uniformly approximated by functions of the form
$\alpha_1+\dots+\alpha_k$, where $\alpha_i\in F(C_i)$ for $L$-curves
$C_1,\dots,C_k$ in some finite covering of $M$. 

Clearly, 5.7(2) implies 5.4(2), and 5.4(1) implies 5.7(1).  It is not
clear if the reverse implications hold, i.e., if Theorems 5.4 and 5.7 
actually express the same dichotomy.

\specialhead 6. Harmonic functions  \endspecialhead

\flushpar Since $\B$ is Stein, every subharmonic function on $X$ extends
to a plurisubharmonic function on $\B$.  Whereas the bounded extension
problem for holomorphic functions is hard to fathom, it is fairly easy
to see, using a little potential theory, that a bounded-above
subharmonic function on $X$ need not extend to a bounded-above
plurisubharmonic function on $\B$. 

Recall that if $\alpha\in\Cal C_\R(\S)$, then there exists a unique
$u\in\Cal C_\R(\overline\B)$ such that
\roster
\item $u$ is plurisubharmonic on $\B$,
\item $(\d\bar\d u)^n=0$, i.e., $u$ is maximal, and
\item $u|\S=\alpha$.
\endroster
Let us write $u=M[\alpha]=M_\B[\alpha]$.  This is the solution of the
Dirichlet problem for the Monge-Amp\gr ere operator, due to Bedford and
Taylor \cite{BT}.  See also earlier work of Bremermann \cite{Bre} and
Walsh \cite{Wal}.  In fact, $u$ is given by the Perron-Bremermann 
formula $u=\sup\Cal F_\alpha$, where $\Cal F_\alpha$ is the set of all 
plurisubharmonic functions $v$ on $\B$ with
$$\limsup_{z\to x} v(z)\leq\alpha(x), \qquad x\in\S.$$
The operator $M:\Cal C_\R(\S)\to\Cal C_\R(\overline\B)$ is continuous in
the sense that if $\alpha_i\to\alpha$ uniformly on $\S$, then $M[\alpha_i]\to
M[\alpha]$ uniformly on $\overline\B$.  Namely, if $\epsilon>0$, then
$\alpha-\epsilon\leq\alpha_i \leq\alpha+\epsilon$ for $i$ sufficiently
large, and then $\Cal F_{\alpha-\epsilon}\subset\Cal F_{\alpha_i}\subset\Cal
F_{\alpha+\epsilon}$, so $M[\alpha]-\epsilon\leq M[\alpha_i]\leq
M[\alpha]+\epsilon$. 

\proclaim{6.1.  Theorem} Let $\alpha\in\Cal C_\R(\S)$.  The following are 
equivalent. 
\roster
\item The harmonic extension $H[\alpha]$ of $\alpha$ to $X$ extends to a
bounded-above plurisubharmonic function on $\B$.
\item $H[\alpha]$ extends to a function in $\Cal C_\R(\overline\B)$ which is
maximal
plurisubharmonic on $\B$.
\item $M[\alpha]$ is harmonic on $X$.
\item $M[\alpha]|X=H[\alpha]$.
\endroster
\endproclaim

\demo{Proof} Clearly, \therosteritem4 $\Leftrightarrow$ \therosteritem3
$\Rightarrow$ \therosteritem2 $\Rightarrow$ \therosteritem1.  We need to
show that \therosteritem1 $\Rightarrow$ \therosteritem4.  Suppose $v$ is
plurisubharmonic and bounded above on $\B$ such that $v|X=H[\alpha]$. 
Now $v$ has a nontangential boundary function $\hat v$, and since
$v|X=H[\alpha]$, we have $\hat v=\widehat{H[\alpha]}=\alpha$ almost
everywhere.  Since $v$ is bounded above and subharmonic, we have
$v=h+p$, where $h$ is harmonic and $p\leq 0$ is a subharmonic potential. 
In fact, $h=H_\B[\alpha]$.  Hence, $$\limsup_{z\to x} v(z) \leq
\limsup_{z\to x} h(z) = \lim_{z\to x} h(z) = \alpha(x)$$ for all
$x\in\S$, so $v\leq M[\alpha]$.  On the other hand, the Perron formulas
for $H[\alpha]$ and $M[\alpha]$ show that $M[\alpha]|X\leq H[\alpha]$. 
Hence, $H[\alpha]=M[\alpha]|X$.  
\qed\enddemo

This shows that if $H[\alpha]$, $\alpha\in\Cal C_\R(\S)$, extends to a
bounded-above plurisubharmonic function on $\B$, then $M[\alpha]$
extends it, so $M[\alpha]|X$ is harmonic.  Clearly, this fails for {\it
most} $\alpha$.  It would be interesting to know if such $\alpha$ can
be the boundary function of the real part of a holomorphic function $f$
on $X$.  Then $f$ would not extend to a bounded holomorphic function on
$\B$. 

\proclaim{6.2.  Corollary} The set of functions $\alpha\in\Cal C_\R(\S)$
such that $M[\alpha]|X$ is harmonic is a closed $\G$-invariant subspace
of $\Cal C_\R(\S)$.  
\endproclaim

\demo{Proof}  Let $\alpha, \beta\in\Cal C_\R(\S)$, and suppose $M[\alpha]$,
$M[\beta]$ are harmonic on $X$.  By the theorem, $H[\alpha]$, $H[\beta]$
extend to bounded-above plurisubharmonic functions on $\B$, but then so
does $H[\alpha+\beta]=H[\alpha]+H[\beta]$, so $M[\alpha+\beta]$ is
harmonic on $X$.  

If $\alpha_i\to\alpha$ uniformly on $\S$, then $M[\alpha_i]\to
M[\alpha]$ uniformly on $\overline\B$, so if $M[\alpha_i]|X$ are harmonic,
then so is $M[\alpha]|X$.  
\qed \enddemo

We obtain a dichotomy analogous to those in Section 5.  If $C$ is an
$L$-curve in a finite covering of $M$ with preimage $X$ in $\B$, let us
denote by $D(C)$ the space of functions $\alpha\in\Cal C_\R(\S)$ such
that $M[\alpha]|X$ is harmonic.  Clearly, $P\subset D(C)$ but, as noted
above, $D(C)$ is considerably smaller than $\Cal C_\R(\S)$. 

\proclaim{6.3. Theorem}  Suppose that $\G$ is arithmetic, and that $L$
is a tensor power of the canonical bundle.  Then one of the 
following holds.
\roster
\item  $D(C)=P$ for every $L$-curve in a finite covering of $M$.
\item  The subspace of $\Cal C_\R(\S)$ generated by $D(C)$ for all
$L$-curves $C$ in finite coverings of $M$ is dense in $\Cal C_\R(\S)$.
\endroster
Suppose furthermore that if $C_1$ and $C_2$ are $L$-curves in the same 
finite covering of $M$, then $D(C_1)=D(C_2)$.  Then $D(C)=P$ for every
$L$-curve $C$ in a finite covering of $M$.
\endproclaim

There are no examples for which it is known which alternative holds in
any of the four dichotomies 5.1, 5.4, 5.7, and 6.3, nor is it known if
these dichotomies are actually different.

\enddocument